\documentstyle{article}

\makeatletter
\def\eqnarray{%
  \stepcounter{equation}%
  \let\@currentlabel=\theequation
  \global\@eqnswtrue
  \global\@eqcnt\z@
  \tabskip\@centering
  \let\\=\@eqncr
  $$\halign to \displaywidth\bgroup\@eqnsel\hskip\@centering
  $\displaystyle\tabskip\z@{##}$&\global\@eqcnt\@ne
  \hfil$\displaystyle{{}##{}}$\hfil
  &\global\@eqcnt\tw@$\displaystyle\tabskip\z@{##}$\hfil
  \tabskip\@centering&\llap{##}\tabskip\z@\cr}
\makeatother

\makeatletter
  \renewcommand{\theequation}{%
        \thesection.\arabic{equation}}
  \@addtoreset{equation}{section}
\makeatother

\begin{document}

\newtheorem{th}{Donotwrite}[section]

\newtheorem{definition}[th]{Definition}
\newtheorem{theorem}[th]{Theorem}
\newtheorem{proposition}[th]{Proposition}
\newtheorem{lemma}[th]{Lemma}
\newtheorem{corollary}[th]{Corollary}
\newtheorem{remark}[th]{Remark}
\newtheorem{example}[th]{Example}

\newfont{\germ}{eufm10}

\def\Aff{\mbox{\sl Aff}\,}
\def\btilde{\tilde{b}}
\def\cd{\cdots}
\def\eps{\epsilon}
\def\goth#1{\mbox{\germ #1}}
\def\ot{\otimes}
\def\Proof{\noindent{\sl Proof.}\quad}
\def\Q{{\bf Q}}
\def\qed{~\rule{1mm}{2.5mm}}
\def\R{{\cal R}}
\def\slchap{\widehat{\goth{sl}}_n}
\def\sln{\goth{sl}_n}
\def\sln-{\goth{sl}_{n-1}}
\def\ub(#1,#2){\underbrace{#1\ot\cd\ot#1}_{#2}}
\def\ubb(#1,#2){\underbrace{#1\cd#1}_{#2}}
\def\veps{\varepsilon}
\def\vphi{\varphi}
\def\wt{\mbox{\sl wt}\,}
\def\xbold{{\bf x}}
\def\Z{{\bf Z}}

\def\vertex(#1,#2,#3,#4){
\setlength{\unitlength}{1mm}
\begin{picture}(20,20)(0,-5)
\put(3,5){\line(1,0){8}}
\put(7,1){\line(0,1){8}}
\put(0,4.3){$#1$}
\put(12,4.3){$#4$}
\put(6.2,10.5){$#2$}
\put(6.2,-2.5){$#3$}
\end{picture}}

\title{ Box Ball System Associated with \\ Antisymmetric Tensor Crystals }

\author{
Daisuke Yamada\thanks{
Department of Mathematical Science,
Graduate School of Engineering Science,
Osaka University, Toyonaka, Osaka 560-8531, Japan}
}

\date{}
\maketitle

\begin{abstract}\noindent
A new box ball system associated with an antisymmetric tensor crystal of the quantum affine algebra
of type A is considered. This includes the so-called colored box ball system with capacity 1 
as the simplest case.
Infinite number of conserved quantities are constructed and the scattering rule of two solitons 
are given explicitly.
\end{abstract}

\section{Introduction}
The box ball system (BBS) \cite{TS,T} is an important example of 
ultra-discrete integrable systems.
It can be obtained from a discrete soliton equation through limiting procedure 
\cite{TTMS,TNS}.
It is also known \cite{HHIKTT,FOY,HKT,HKOTY} that the BBS has 
an equivalent formulation in terms of affine crystal theory \cite{K,KMN}. 
Following \cite{FOY} we briefly review the BBS in this formulation.
For $l\in\Z_{\ge1}$ set 
\[
B_l=\{[i_1i_2\cd i_l]\mid i_k\in\Z,1\le i_1\le\cd\le i_l\le n\}.
\]
Here $n$ is a fixed integer greater than 1. 
$B_l$ is identified with the crystal of the symmetric tensor representation of 
degree $l$ of the quantum affine algebra $U'_q(\widehat{sl}_n)$.
By definition, the combinatorial $R$ for the tensor product crystal 
$B_{l_1}\ot B_{l_2}$ is a bijection
\[
R:\quad B_{l_1}\ot B_{l_2}\longrightarrow B_{l_2}\ot B_{l_1}
\]
which commutes with the crystal operators \cite{KMN}. 
Use the symbol 
\begin{center}
\vertex(b_1,b_2,\tilde{b}_2,\tilde{b}_1)
\end{center}
if $R(b_1\ot b_2)=\tilde{b}_2\ot\tilde{b}_1$. 
Applying the combinatorial $R:\,B_3\ot B_1\rightarrow B_1\ot B_3$ successively, we have the following diagram.

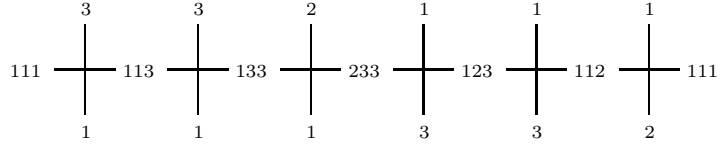
\setlength{\unitlength}{1mm}
\begin{picture}(120,30)(0,0)
\multiput(11,15)(15,0){6}{\line(1,0){8}}
\multiput(15,9)(15,0){6}{\line(0,1){12}}

\put(5,14){$\scriptstyle 111$}
\put(20,14){$\scriptstyle 113$}
\put(35,14){$\scriptstyle 133$}
\put(50,14){$\scriptstyle 233$}
\put(65,14){$\scriptstyle 123$}
\put(80,14){$\scriptstyle 112$}
\put(95,14){$\scriptstyle 111$}

\put(14.4,22.3){$\scriptstyle 3$}
\put(29.4,22.3){$\scriptstyle 3$}
\put(44.4,22.3){$\scriptstyle 2$}
\put(59.4,22.3){$\scriptstyle 1$}
\put(74.4,22.3){$\scriptstyle 1$}
\put(89.4,22.3){$\scriptstyle 1$}

\put(14.4,5.9){$\scriptstyle 1$}
\put(29.4,5.9){$\scriptstyle 1$}
\put(44.4,5.9){$\scriptstyle 1$}
\put(59.4,5.9){$\scriptstyle 3$}
\put(74.4,5.9){$\scriptstyle 3$}
\put(89.4,5.9){$\scriptstyle 2$}
\end{picture}

\noindent
Here we omitted the symbol $[\;]$. 
For instance, $113$ (resp. $2$) means $[113]$ (resp. $[2]$).
Neglecting the horizontal line, 
we see the upper state is mapped to the lower one as the time is evolved. 
Each $[\cdot]$ corresponds to a box. 
We consider $[1]$ as a vacant box, and $[2],[3],\cd$ as boxes with a ball.
The above example shows that the array $[3][3][2]$ behaves as a soliton 
just like in the classical soliton theory. 
One can also check that the longer soliton moves faster. 
Therefore, if there are two solitons of different lengths, 
we can expect a scattering of solitons. 
The following results are shown in \cite{FOY}. (See also \cite{TNS,HHIKTT}.)

\begin{itemize}
\item[(1)] For an element $[i_1i_2\cd i_l]$ of $B_l$ such that $i_1\ge2$, 
one can associate a soliton state $[i_l]\cd[i_2][i_1]$.
\item[(2)] The scattering of two solitons of lengths $l_1$ and 
$l_2$ ($l_1>l_2$) is described by the combinatorial 
$R:\,B_{l_1}\ot B_{l_2}\rightarrow B_{l_2}\ot B_{l_1}$.
\item[(3)] The phase shift caused by the scattering is described by 
the energy function, an integer valued function canonically defined on 
the tensor product crystal $B_{l_1}\ot B_{l_2}$.
\end{itemize}

We want to extend these results to more general cases.
It is known that the set of semi-standard tableaux $B^{k,l}$ of shape $(l^k)$ 
with letters in $\{1,2,\cd,n\}$ also admits the affine crystal structure. 
Considering the combinatorial 
$R:\,B^{2,3}\ot B^{2,1}\rightarrow B^{2,1}\ot B^{2,3}$, 
we obtain the following.

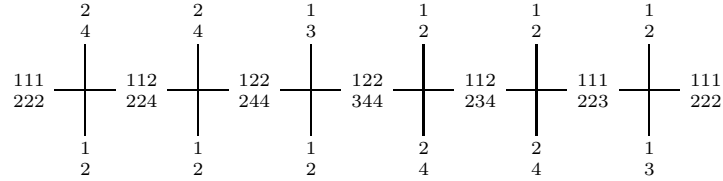
\setlength{\unitlength}{1mm}
\begin{picture}(120,30)(0,0)
\multiput(11,15)(15,0){6}{\line(1,0){8}}
\multiput(15,9)(15,0){6}{\line(0,1){12}}

\put(5,14){${111\atop222}$}
\put(20,14){${112\atop224}$}
\put(35,14){${122\atop244}$}
\put(50,14){${122\atop344}$}
\put(65,14){${112\atop234}$}
\put(80,14){${111\atop223}$}
\put(95,14){${111\atop222}$}

\put(13.9,23.3){${2\atop4}$}
\put(28.9,23.3){${2\atop4}$}
\put(43.9,23.3){${1\atop3}$}
\put(58.9,23.3){${1\atop2}$}
\put(73.9,23.3){${1\atop2}$}
\put(88.9,23.3){${1\atop2}$}

\put(13.9,4.9){${1\atop2}$}
\put(28.9,4.9){${1\atop2}$}
\put(43.9,4.9){${1\atop2}$}
\put(58.9,4.9){${2\atop4}$}
\put(73.9,4.9){${2\atop4}$}
\put(88.9,4.9){${1\atop3}$}
\end{picture}

\noindent
As seen in this example, 
the situation for general $k$ is similar to the previous case 
corresponding to $k=1$. 
In this paper, we are to show the following results.

\begin{itemize}
\item[(1)] One can associate a soliton state for an element of 
$B^{k-1,l}$(letters in $\{1,2,\cd,k\}$)
$\times B^{1,l}$(letters in $\{k+1,k+2,\cd,n\}$).
\item[(2)] The scattering of two solitons of lengths 
$l_1$ and $l_2$ ($l_1>l_2$) is described by the product of combinatorial $R$'s
\[
R\times R:\;(B^{k-1,l_1}\ot B^{k-1,l_2})\times(B^{1,l_1}\ot B^{1,l_2})
\longrightarrow (B^{k-1,l_2}\ot B^{k-1,l_1})\times(B^{1,l_2}\ot B^{1,l_1}).
\]
\item[(3)] The phase shift is described by the sum of the energy functions 
corresponding to the tensor product crystals $B^{k-1,l_1}\ot B^{k-1,l_2}$ and 
$B^{1,l_1}\ot B^{1,l_2}$.
\end{itemize}

\section{Box Ball System}
In this section we review the minimum contents concerned with the 
$U'_q (\widehat{sl}_n)$-crystal $B^{k,l}$ 
and define the box ball system we shall investigate.

\subsection{The crystal $B^{k,l}$} \label{subsec:crystal}
Fix an integer $n\in \Z_{\ge2}$. 
For arbitrary integers $k,l$ such that $1\le k\le n-1, l\ge1$, 
let $B^{k,l}$ be the set of semi-standard tableaux of shape $(l^k)$ 
with letters in $\{1,2,\ldots,n\}$. 
A tableau is described by using ${\;\brack\;}$ and ${}^t{\;\brack\;}$
signifies its transpose. 
For instance,
\[
\left[\begin{array}{c} 1 \\ 3 \\ 5 \end{array}\right]\in B^{3,1},\quad
\left[\begin{array}{ccc} 1&1&2 \\ 2&3&3 \end{array}\right]\in B^{2,3}\quad\mbox{and}\quad
{}^t\left[\begin{array}{ccc} 1&1&2 \\ 2&3&3 \end{array}\right]=
\left[\begin{array}{cc} 1&2 \\ 1&3 \\ 2&3 \end{array}\right]\in B^{3,2}.
\]
It is known \cite{KMN2,S} that $B^{k,l}$ admits the structure of 
$U'_q(\widehat{sl}_n)$-crystals, i.e., one has the action of operators $e_i,f_i$, called 
Kashiwara operators, 
\[
e_i,f_i :B^{k,l}\longrightarrow B^{k,l}\sqcup \{ 0\} \quad\mbox{for }i=0,1,\ldots,n-1.
\]
We are to explain the actions of $e_i,f_i$ on $B^{k,l}$ for $i\neq0$ in detail. To do this we 
first give the actions on $(B^{k,1})^{\otimes L}$, $L$-th tensor power of $B^{k,1}$.
The rule is given as follows.
\begin{itemize}
\item[(1)] We identify $b=\ot_{j=1}^L b_j\in (B^{k,1})^{\ot L}$ with $b'\in(B^{1,1})^{\ot kL}$ 
in such a way that
\[
b'=\ot _{j=1}^{L} (\ot _{m=1}^{k} [x_m^j])\quad\mbox{if }b_j={}^t[x^j_1,x^j_2,\ldots,x^j_k].
\]
Here and in what follows $\ot_{j=1}^L b_j$ means $b_1\ot b_2\ot\cd\ot b_L$.
\item[(2)] Reading the letters $x^j_m$ in $b'$ from left to right, we construct a sequence of
$\pm$ or $0$ by associating $+$ if $x^j_m=i$, $-$ if $x^j_m=i+1$, and $0$ otherwise.
\item[(3)] Neglecting all $0$, replace adjacent $+-$ pairs with $00$ successively until we get 
the following sequence.
\[
--- \cdots --+++ \cdots ++
\]
\item[(4)] The action of $e_i$ (resp. $f_i$) on $b$ amounts to replacing the rightmost $-$
(resp. leftmost $+$) of the above sequence with $+$ (resp. $-$), i.e., changing the letter
$i+1$ (resp. $i$) to $i$ (resp. $i+1$) on the corresponding tensor component of $b'$.
\end{itemize}

\begin{example}
\begin{eqnarray*}
f_2 (
\left[
\begin{array}{c}
1\\
2
\end{array}
\right]
\ot 
\left[
\begin{array}{c}
2\\
3
\end{array}
\right]
\ot 
\left[
\begin{array}{c}
4\\
6
\end{array}
\right]
) 
&=& f_2 (1 2 2 3 4 6) 
= f_2 (0 + + - 0 0) 
= f_2 (0 + 0 0 0 0) \\
&=& 0 - 0 0 0 0 
=\left[
\begin{array}{c}
1\\
3
\end{array}
\right]
\ot 
\left[
\begin{array}{c}
2\\
3
\end{array}
\right]
\ot 
\left[
\begin{array}{c}
4\\
6
\end{array}
\right]
\end{eqnarray*}
\end{example} 

To compute the action on $B^{k,l}$ we consider the map 
\begin{eqnarray*}
\mbox{sp}:B^{k,l}&\longrightarrow&(B^{k,1})^{\ot l}\\
b=[x^1,x^2,\cd,x^l]&\mapsto&x^l\ot\cd\ot x^2\ot x^1,
\end{eqnarray*}
where $x^j$ stands for the $j$-th column of $b$. Then the action of $e_i$ (resp. $f_i$)
on $B^{k,l}$ is given by $\mbox{sp}^{-1}\circ e_i\circ\mbox{sp}$ (resp. 
$\mbox{sp}^{-1}\circ f_i\circ\mbox{sp}$).

\begin{remark}
For the action of $e_0,f_0$ we refer to \cite{S}, since it is unnecessary for our purpose.
\end{remark}

This rule, called {\rm signature rule}, can be applicable to any tensor product of the form
$B^{k_1,l_1}\ot\cd\ot B^{k_d,l_d}$ by embedding it into 
$(B^{k_1,1})^{\ot l_1}\ot\cd\ot(B^{k_d,1})^{\ot l_d}$ using $\mbox{sp}^{\ot d}$.
Although it is inhomogeneous in the sense that $k_j$ can vary, 
the signature rule holds as it is.

\subsection{Combinatorial $R$ and energy function}
The combinatorial $R$ for the tensor product crystal $B_1\ot B_2$ is, 
by definition, a bijection
$R:B_1\ot B_2\longrightarrow B_2\ot B_1$
which commutes with the operators $e_i$ and $f_i$ for all $i$. 
For our case concerning the $U'_q(\widehat{sl}_n)$-crystal $B^{k,l}$, 
such $R$ is known to exist uniquely.
To describe the combinatorial $R$ explicitly, 
we explain Schensted's bumping algorithm. See e.g. \cite{F} for the details. 
For 
\[
x=
\left[
\begin{array}{c}
x_1    \\
x_2    \\
\vdots \\
x_k
\end{array}
\right]
=\left[
\begin{array}{cccc}
x_1^1  & x_1^2  & \cd    & x_1^l  \\
x_2^1  & x_2^2  & \cd    & x_2^l  \\
\vdots & \vdots & \ddots & \vdots \\
x_k^1  & x_k^2  & \cd    & x_k^l 
\end{array}
\right]
\in B^{k,l},
\]
we define a row word $\mathrm{row}(x)$ by 
\[
\mathrm{row}(x)=
\underbrace{x_k^1 \ldots x_k^l}_{x_k} \ldots  \cdots \ldots 
\underbrace{x_2^1 \ldots x_2^l}_{x_2} 
\underbrace{x_1^1 \ldots x_1^l}_{x_1}.
\]

Now we recall the bumping algorithm (row-bumping, or row-insertion), 
for constructing a new tableau from a tableau by inserting an integer. 
For inserting an integer $i$ in a tableau $T$, the rule of bumping 
"$T\leftarrow i$" is given as follows. 
\begin{itemize}
\item[(1)] If there are no integers larger than $i$ in the first row, 
add a new empty box at the right end, and put $i$ in it. 
\item[(2)] Otherwize, among the integers larger than $i$, 
find the leftmost one, say $j$, and replace $j$ with $i$. 
Then inset $j$ into the second row in the same way. 
\item[(3)] Repeat this procedure until the bumped number can be put in a new 
box at the right end of the row.  
\end{itemize}

The following is the explicite algorithm of combinatrial $R$ on 
$B^{k,l}\ot B^{k',l'}$.  
\begin{theorem}(\cite{S}) \label{th:comb R}
$x\ot y$ is mapped to $\tilde{x}\ot\tilde{y}$ by the combinatorial $R$
\[
B^{k,l}\ot B^{k',l'}\longrightarrow B^{k',l'}\ot B^{k,l},
\]
if and only if
\[
y \leftarrow \mathrm{row}(x) =\tilde{y} \leftarrow \mathrm{row}(\tilde{x}).
\]
\end{theorem}

Next we define the energy function $H$ on $B^{k,l}\ot B^{k',l'}$ as follows. 
\begin{definition} \label{def:H}
Let $x\in B^{k,l},y\in B^{k',l'}$. Let $d(x,y)$ be the number of nodes 
in the shape of 
$y\leftarrow \mathrm{row}(x)$ that are strictly east of the $\max(l,l')$-th column. 
Then the energy
function $H:B^{k,l}\ot B^{k',l'}\rightarrow\Z$ is given by 
\[
H(x\ot y)=d(x,y)-\min(k,k')\min(l,l').
\]
\end{definition}

\begin{remark}
Intrinsically, the energy function is uniquely defined by specifying the difference
$H(e_i(x\ot y))-H(x\ot y)$ up to an additive constant. See section 4.1 of \cite{KMN}.
The above rule for the energy function is established in \cite{S,SW}. 
The normalization
of $H$ in Definition \ref{def:H} is so fixed that the maximum is 0.
\end{remark}

Here we give an example for combinatrial $R$ on $B^{3,3} \ot B^{2,1}$ 
and energy function assosiate with this $R$. 
\begin{example}
$$
\left[
\begin{array}{ccc}
1 & 2 & 4 \\
2 & 3 & 5 \\
4 & 4 & 6
\end{array}
\right] 
\ot 
\left[
\begin{array}{c}
2 \\
5
\end{array}
\right] 
\simeq 
\left[
\begin{array}{c}
2 \\
4
\end{array}
\right] 
\ot 
\left[
\begin{array}{ccc}
1 & 2 & 3 \\
2 & 4 & 5 \\
4 & 5 & 6
\end{array}
\right] 
$$
\end{example}

We will show the procedure for culculating the right hand side 
when left hand side is given. We set 
\[
x=\left[
\begin{array}{ccc}
1 & 2 & 4 \\
2 & 3 & 5 \\
4 & 4 & 6
\end{array}
\right], \quad 
y=\left[
\begin{array}{c}
2 \\
5
\end{array}
\right].  
\]
By the rule of bumping we have the following tableau. 
\[
y\leftarrow \mathrm{row}(x)=
\begin{array}{cc}
2 & \leftarrow 446235124 \\
5 & {}
\end{array}
=
\begin{array}{cccc}
1 &  2 &  2 &  4 \\
2 &  3 &  5 & {} \\
4 &  4 &  6 & {} \\
5 & {} & {} & {}
\end{array} 
\]
Note that the rule of bumping is reversible procedure, and we can culculate 
${\tilde x}\in B^{2,1}$ and ${\tilde y}\in B^{3,3}$ such that 
$x\ot y\simeq {\tilde x}\ot {\tilde y}$. 

\begin{eqnarray*}
\begin{array}{cccc}
1 &  2 &  2 &  4 \\
2 &  3 &  5 & {} \\
4 &  4 &  6 & {} \\
5 & {} & {} & {}
\end{array} 
&=&
\begin{array}{ccccc}
1 &  2 &  2 &          4 & {} \\
2 &  4 &  5 &         {} & {} \\
4 &  4 &  6 & \leftarrow & 4 
\end{array}
=
\begin{array}{ccccc}
1 &  2 &  2 &          4 & {} \\
2 &  4 &  5 & \leftarrow &  3 \\
4 &  4 &  6 &         {} & {} 
\end{array}
=
\begin{array}{cccccc}
1 &  2 &  3 &  4 & \leftarrow &  2 \\
2 &  4 &  5 & {} &         {} & {} \\
4 &  4 &  6 & {} &         {} & {}
\end{array}
\\
&=&
\begin{array}{ccccc}
1 &  2 &  3 & \leftarrow & 42 \\
2 &  4 &  5 &        {}  & {} \\
4 &  4 &  6 &        {}  & {}
\end{array}
={\tilde y}\leftarrow \mathrm{row} ({\tilde x})
\end{eqnarray*}

According to the definition of energy function, we have 
\[
H(x\ot y)=1-\mathrm{min}(3,2) \mathrm{min}(3,1) =-1. 
\]

\subsection{Box ball system}
In this subsection we fix $k$ ($1\le k\le n-1$) and set $B_l =B^{k,l}$. We employ the
following notations.

\[
1_k = \left[
\begin{array}{c}
1\\ 2\\ \vdots\\ k
\end{array}\right]\in B_1,\quad
c_l = 1_k^l = \left[
\begin{array}{cccc}
1 & 1 & \cd & 1 \\ 
2 & 2 & \cd & 2 \\ 
\vdots & \vdots & \ddots & \vdots\\ 
k & k & \cd & k
\end{array}\right]
\in B_l.
\]
Consider the crystal $B_1^{\ot L}$ with $L$ sufficiently large. We call an element $p$ of 
$B_1^{\ot L}$ a state, if
\[
p=(\ot_{j=1}^d b_j)\ot 1_k^{\ot(L-d)}
\]
and $L$ is sufficiently large compared with $d$. Note that we allow $1_k$ to be in 
$\ot_{j=1}^d b_j$.

\begin{lemma} \label{lem:stab}
By iterating the combinatorial $R:\,B_l\ot B_1\rightarrow B_1\ot B_l$ L times, we have
the map 
\[
B_l\ot B_1^{\ot L} \longrightarrow B_1^{\ot L} \ot B_l.
\]
Suppose $c_l \ot p$ is mapped to $\tilde{p} \ot c$ by this map. Then, 
if $L\gg 1$, we have $c=c_l$. 
\end{lemma}
By this lemma one can define the time evolution operator $T_l$ by
\[
T_l :B_1^{\ot L} \ni p \mapsto \tilde{p} \in B_1^{\ot L}.
\]
Then we have
\begin{lemma} \label{lem:T_l commutes}
$[T_l, e_i]=[T_l, f_i]=0$ \quad for $i \in I\setminus\{0,k\}$.
\end{lemma}
Proof. By the definition of $T_l$, we have 
\[
c_l \ot e_i p \mapsto T_l (e_i p)\ot c_l, \quad
e_i (c_l \ot p) \mapsto e_i(T_l (p)\ot c_l),
\]
under the map $B_l\ot B_1^{\ot L}\rightarrow B_1^{\ot L}\ot B_l$. 
On the other hand, for $i\in I\setminus\{0,k\}$ using the signature rule 
in section \ref{subsec:crystal} we have 
\[
e_i (c_l \ot p)= c_l \ot e_i p ,\quad e_i(T_l (p)\ot c_l)= e_i T_l (p)\ot c_l. 
\]
Thus we have $T_l (e_i p) \ot c_l = e_i T_l (p) \ot c_l$. 
The case for $f_i$ is similar.
$\Box$
{}This property will be used essentially to prove our main theorem.

We also define a map $E_l :B_1^{\ot L}\rightarrow {\bf Z}_{\ge0}$ by
\begin{equation} \label{eq:def_E}
E_l(p)=-\sum_{j=1}^{L}H(b^{(j-1)} \ot b_j),
\end{equation}
where $b^{(j)}$ $(0\leq j<L)$ is defined by
\[
B_l \ot B_1^{\ot j} \ni c_l \ot (\ot_{i=1}^{j} b_i) 
\mapsto 
(\ot_{i=1}^{j} \tilde{b}_i)\ot b^{(j)} \in B_1^{\ot j} \ot B_l
\]
and $H$ is the energy function on $B_l\ot B_1$. 
Note that for a state $p$, $E_l(p)$ does not depend on $L$, 
because of Lemma \ref{lem:stab} and the normalization $H(1_k^l\ot 1_k)=0$.

It is known that the following equation, called the Yang-Baxter equation, 
holds on $B_l \ot B_{l'} \ot B_{l''}$.
\[
(R\ot 1)(1\ot R)(R\ot 1)=(1\ot R)(R\ot 1)(1\ot R)
\]
{}From this fact one can show

\begin{theorem}(\cite{FOY})
\begin{itemize}
\item[(1)] The commutation relation among the time evolutions:\quad  $T_l T_{l'} (p)=T_{l'} T_l(p)$.
\item[(2)] The energy conservation:\quad  $E_l (T_{l'}(p))=E_l (p)$.
\end{itemize}
\end{theorem}
To be precise, we need the Yang-Baxter equation for the affinization of $B_l,B_{l'},B_{l''}$
to prove (2). See section 2.4 of \cite{FOY}.

In what follows in this section, we present other conserved quantities $P_{\le k},P_{>k}$
under time evolutions constructed by Schensted's bumping algorithm. By the isomorphism 
$B_l \ot B_1^{\ot L} \simeq B_1^{\ot L} \ot B_l$, we have  
\[
c_l \ot (\ot_{j=1}^L b_j)\mapsto (\ot_{j=1}^L {\tilde b_j})\ot c_l.
\]
In view of Theorem \ref{th:comb R} one notices that the tableau constructed by
\[
((\cd(\mathrm{row}(b_L)\leftarrow \mathrm{row}(b_{L-1}))\leftarrow \cdots )\leftarrow \mathrm{row}(b_1))\leftarrow \mathrm{row}(c_l)
\]
and the one by
\[
((\cd(\mathrm{row}(c_l)\leftarrow \mathrm{row}(\tilde{b}_L))\cd)\leftarrow \mathrm{row}(\tilde{b}_2))\leftarrow \mathrm{row}(\tilde{b}_1)
\]
are the same. It means that the corresponding words are Knuth-equivalent \cite{F}.
\begin{equation} \label{eq:row words}
\mathrm{row}(b_L)\,\mathrm{row}(b_{L-1})\,\cd\,\mathrm{row}(b_1)\,\mathrm{row}(c_l)\sim
\mathrm{row}(c_l)\,\mathrm{row}(\tilde{b}_L)\,\cd\,\mathrm{row}(\tilde{b}_2)\,\mathrm{row}(\tilde{b}_1).
\end{equation}
Remove the numbers $k+1,k+2,\ldots,n$ from both sequences, 
then they are still Knuth-equivalent. (See e.g. Lemma 1 of \cite{Fu}.) Thus the tableaux constructed 
by both words are the same. It gives a conserved quantity, denoted by $P_{\le k}$, under time evolutions.
Note that taking $\mathrm{row}(c_l)$ in the words corresponds to adding $c_l$ from the left at the
construction of tableaux. Similarly, by removing the numbers $1,2,\ldots,k$ from (\ref{eq:row words}),
we obtain another conserved quantity $P_{>k}$.
Note that the $P$-tableau in Theorem 3.1 of \cite{Fu} corresponds to $P_{>k}$.

\section{Solitons}
In the previous section we introduced a box ball system on the crystal $(B^{k,1})^{\ot L}$. 
In this section we define soliton states and investigate the scattering of two 
solitons.

\subsection{One soliton state}
Consider the following state 
\[
p=1_k^{\ot c}\ot (\ot_{j=1}^d b_j)\ot 1_k^{\ot (L-c-d)}\in (B^{k,1})^{\ot L}, 
\quad 
b_j={}^t[x_1^j,x_2^j,\ldots,x_k^j]\in B^{k,1}.
\]

\begin{definition} \label{def:one soliton}
We call such $p$ a one soliton state, if the integers 
$x_i^j\in\{1,2,\ldots,n\}$ ($1\le i\le k,1\le j\le d$) satisfy 
\begin{eqnarray*}
&&x_i^j\ge x_i^{j+1}\quad(1\le i\le k,1\le j\le d-1),\\
&&x_i^j<x_{i+1}^j\quad(1\le i\le k-2,1\le j\le d),\\
&&x_{k-1}^j\le k<x_k^j\quad(1\le j\le d).
\end{eqnarray*}
\end{definition}
We call $d$ the length of the soliton.

\begin{remark} \label{rem:E_1=1}
Such state is characterized by the condition $E_1(p)=1$.
\end{remark}

\begin{example}
In the case of $k=3$, the following gives a one soliton state of length 3.
\[
\cdots\ot
\left[
\begin{array}{c}
1\\
2\\
3
\end{array}
\right]
\ot
\left[
\begin{array}{c}
2\\
3\\
5
\end{array}
\right]
\ot
\left[
\begin{array}{c}
2\\
3\\
4
\end{array}
\right]
\ot
\left[
\begin{array}{c}
1\\
2\\
4
\end{array}
\right]
\ot
\left[
\begin{array}{c}
1\\
2\\
3
\end{array}
\right]
\ot\cdots
\]
\end{example}

To simplify the notation, we use
\begin{eqnarray*}
&&[{\bf 0}]={}^t[1,\ldots,k-2,k-1,k],\\
&&[{\bf 1}]={}^t[1,\ldots,k-2,k-1,k+1],\\
&&[{\bf 2_{\pm}}]={}^t[1,\ldots,k-2,k-\frac12\pm\frac12,k+\frac32\mp\frac12],\\
&&[{\bf 3}]={}^t[1,\ldots,k-2,k,k+2],\\
&&[{\bf 4}]={}^t[1,\ldots,k-2,k+1,k+2].
\end{eqnarray*}
The following lemma shows how one soliton state behaves as the time is evolved.

\begin{lemma}
Let $p=1_k^{\ot c}\ot S_d\ot 1_k^{\ot (L-c-d)}$ be a one soliton state of length $d$. Then,
\[
T_l^t(p)=1_k^{\ot (c+\min (d,l)t)}\ot S_d\ot1_k^{\ot (L-c-d-\min (d,l)t)}.
\]
\end{lemma}
Proof. 
By Lemma \ref{lem:T_l commutes}, it is sufficient to check for the elements
killed by $e_i$ for all $i\in I\setminus\{0,k\}$, i.e..
\[
p=[{\bf 0}]^{\ot c}\ot[{\bf 1}]^{\ot d}\ot[{\bf 0}]^{\ot (L-c-d)}.
\]
Our $R:\,B_l\ot B_1\rightarrow B_1\ot B_l$ satisfies
\begin{eqnarray*}
&&R(\xi_i\ot[{\bf 1}])=[{\bf 0}]\ot\xi_{i+1}\;(0\le i<l),\quad
R(\xi_l\ot[{\bf 1}])=[{\bf 1}]\ot\xi_l,\\
&&R(\xi_i\ot[{\bf 0}])=[{\bf 1}]\ot\xi_{i-1}\;(0<i\le l),\quad
R(\xi_0\ot[{\bf 0}])=[{\bf 0}]\ot\xi_0,
\end{eqnarray*}
where $\xi_i=[{\bf 0}^{l-i} {\bf 1}^i]\quad (0\leq i\leq l)$. The statement is clear from these formulas.
$\Box$

We call the value $c$ the phase of the soliton $S_d$. Note that the phase is related to the position
of the soliton {\em before} applying time evolutions. Hence it does not depend on time for a one
soliton state.

\subsection{The scattering rule of solitons}
Imitating the definition of one soliton state, we can consider an $m$ soliton state 
\begin{equation} \label{eq:m soliton}
p_{ms} =\ot_{j=1}^m (1_k^{\ot(c_j-c_{j-1}-d_{j-1})}\ot S_{d_j})\ot 1_k^{\ot(L-c_m-d_m)},
\end{equation}
where $c_0=d_0=0$. $S_{d_j}$ signifies a soliton of length $d_j$ satisfying the conditions in
Definition \ref{def:one soliton}, and we assume $c_j-c_{j-1}-d_{j-1}$ ($1\leq j\leq m$) are
relatively large compared with $d_j$ ($j=1,\ldots,m$).

For a single soliton we introduce a new notation.
One is related to the internal degree of freedom.
For $S_d =\ot _{j=1}^d b_j$  we set $u=[b_d,\ldots,b_2,b_1]$.
Note that from Definition \ref{def:one soliton} $u$ can be regarded as a tableau in $B^{k,l}$.
Next we split $u$ between the $(k-1)$-th and $k$-th row as 
\begin{equation} \label{eq:u split}
u=(u_{< k},u_k)\in B^{k-1,l} \times B^{1,l}.
\end{equation}
Namely $u_k$ is the $k$-th row of $u$ and $u_{< k}$ is the tableau of $k-1$ rows obtained by
removing the $k$-th row of $u$. Note that $u_{<k}$ (resp. $u_k$) is a tableau with letters in
$\{1,2,\ldots,k\}$ (resp. $\{k+1,k+2,\ldots,n\}$). Thus $B^{k-1,l}$ (resp. $B^{1,l}$) in 
(\ref{eq:u split}) should be considered as a $U_q(sl_k)(=<e_i,f_i,t_i
(i=1,2,\ldots,k-1)>)$-crystal (resp. $U_q(sl_{n-k})(=<e_i,f_i,t_i(i=k+1,k+2,\ldots,n-1)>)$-crystal).
The other is related to the phase explained in the end of the previous subsection. If the soliton
with the internal degree of freedom $u$ has a phase $c$, we represent it as $\zeta^c u$ with 
$\zeta$ be an indeterminate.

Next we consider a 2 soliton state. Let $\zeta^{c_1}u$ and $\zeta^{c_2}v$ be solitons with
notations above. We assume the length $d_1$ of the soliton $\zeta^{c_1}u$ is bigger than
the length $d_2$ of $\zeta^{c_2}v$. We also assume that the former is situated far left to the latter
at time $t=0$. We represent such situation as $\zeta^{c_1}u \ot \zeta^{c_2}v$. Now apply $T_l$ 
($l>d_2$) many times. Then one observes that the longer soliton passes the shorter one
and the phases change by constant. We call this phenomenon a scattering of 2 solitons and represent as
\[
\zeta^{c_1} u\ot \zeta^{c_2} v \longrightarrow 
\zeta^{c_2-{\delta}_2} {\tilde v}\ot \zeta^{c_1+{\delta}_1}{\tilde u}.
\]

Then the following gives the scattering rule of such 2 solitons.

\begin{theorem}
The scattering rule of 2 soliton state $p_{2s}(d_1>d_2)$ 
under the time evolution $T_l\;(l>d_2)$ is described as follows.
The change of the internal degree of freedom is given by the product of 
combinatorial $R$'s $R\times R$:
\begin{eqnarray*}
(B^{k-1,d_1} \ot B^{k-1,d_2}) \times (B^{1,d_1} \ot B^{1,d_2})
&\longrightarrow& 
(B^{k-1,d_2} \ot B^{k-1,d_1}) \times (B^{1,d_2} \ot B^{1,d_1}) \\
(u_{< k} \ot v_{< k},u_k \ot v_k) 
&\mapsto& 
({\tilde v}_{< k}\ot {\tilde u}_{< k},{\tilde v}_k \ot {\tilde u}_k)
\end{eqnarray*}
The phase shifts $\delta_1,\delta_2$ are given by the sum of 
the energy functions $H$ associated with each combinatorial $R$:
\[
{\delta}_1 ={\delta}_2 =2d_2 +H(u_k \ot v_k)+H(u_{< k} \ot v_{< k})
\]
\end{theorem}
The proof will be given in the next subsection.

\begin{remark}
For a 2 soliton state $p_{2s}$ we regard $c_2 -c_1 -d_1$ as the distance from $S_{d_1}$ to $S_{d_2}$.
The above scattering rule of 2 solitons is valid if $c_2 -c_1 -d_1 \geq d_2$ 
holds at the initial stage.
\end{remark}

\begin{remark}
The image $(R\times R)(u_{<k}\ot v_{<k},u_k\ot v_k)$ is different from $R(u\ot v)$ by 
$R:B^{k,d_1}\ot B^{k,d_2}\rightarrow B^{k,d_2}\ot B^{k,d_1}$. For instance, we have
\begin{eqnarray*}
(R\times R)&(&
\left[\begin{array}{ccccc}
1&1&1&1&2\\
2&2&3&3&3
\end{array}\right]
\ot
\left[\begin{array}{ccc}
1&1&2\\
2&3&3
\end{array}\right],
\left[\begin{array}{ccccc}
4&4&4&5&5
\end{array}\right]
\ot
\left[\begin{array}{ccc}
5&6&7
\end{array}\right])\\
&=&
(\left[\begin{array}{ccc}
1&1&2\\
2&3&3
\end{array}\right]
\ot
\left[\begin{array}{ccccc}
1&1&1&1&2\\
2&2&3&3&3
\end{array}\right],
\left[\begin{array}{ccc}
4&5&5
\end{array}\right]
\ot
\left[\begin{array}{ccccc}
4&4&5&6&7
\end{array}\right]),
\end{eqnarray*}
whereas
\[
R(\left[\begin{array}{ccccc}
1&1&1&1&2\\
2&2&3&3&3\\
4&4&4&5&5
\end{array}\right]
\ot
\left[\begin{array}{ccc}
1&1&2\\
2&3&3\\
5&6&7
\end{array}\right])
=
\left[\begin{array}{ccc}
1&1&2\\
3&3&3\\
4&5&5
\end{array}\right]
\ot
\left[\begin{array}{ccccc}
1&1&1&1&2\\
2&2&2&4&4\\
3&3&5&6&7
\end{array}\right].
\]
\end{remark}

Recall the scattering rule of 2 solitons is written as
\[
\zeta^{c_1} 
\left[
\begin{array}{c}
u_{< k} \\
u_k
\end{array}
\right]
\ot \zeta^{c_2}
\left[
\begin{array}{c}
v_{< k} \\
v_k
\end{array}
\right]
\rightarrow \zeta^{c_2 -\delta}
\left[
\begin{array}{c}
\tilde{v}_{< k} \\
\tilde{v}_k
\end{array}
\right]
\ot \zeta^{c_1 +\delta}
\left[
\begin{array}{c}
\tilde{u}_{< k} \\
\tilde{u}_k
\end{array}
\right]
\]
\[
\delta =2d_2 +H(u_k \ot v_k )+H(u_{< k} \ot v_{< k}).@
\]
Denote this map by ${\tilde R}$. Then it satisfies the Yang-Baxter equation

\begin{equation} \label{eq:YBtilde}
(\tilde{R} \ot 1)(1\ot \tilde{R})(\tilde{R} \ot 1)=(1\ot \tilde{R})(\tilde{R} \ot 1)(1\ot \tilde{R}).
\end{equation}
It means that the scattering rule of 3 solitons is independent of the order of 
2 body ones. 
It is so even when 3 solitons collide almost at the same time.
We check this property with the following example. 
Here we employ the notation
\[
(x_1,x_2,\ldots,x_n)=\ot _{j=1}^{n} x_j \ot 1_k \ot 1_k \ot\cdots.
\]

\begin{example}
\begin{eqnarray*}
p_{3s}&=&
\left(
\begin{array}{cccccccccccc}
2&2&2&1&1&1&2&1&1&1&1&1\\
3&3&3&2&2&2&3&2&2&2&2&3\\
5&4&4&3&3&3&6&4&3&3&3&5
\end{array}
\right)
\\
T_3^1 (p_{3s})&=&1_3^{\ot 3}\ot
\left(
\begin{array}{ccccccccccc}
2&2&2&1&1&2&1&1&1&1\\
3&3&3&2&2&3&2&2&2&3\\
5&4&4&3&3&6&4&3&3&5
\end{array}
\right)
\\
T_3^2 (p_{3s})&=&1_3^{\ot 6}\ot
\left(
\begin{array}{cccccccc}
2&2&2&1&2&1&1&1\\
3&3&3&2&3&2&2&3\\
5&4&4&3&6&4&3&5
\end{array}
\right)
\\
T_3^3(p_{3s})&=&1_3^{\ot 9}\ot
\left(
\begin{array}{cccccc}
2&2&1&2&2&1\\
3&3&2&3&3&4\\
5&4&3&6&4&5
\end{array}
\right)
\\
T_3^4 (p_{3s})&=&1_3^{\ot 11}\ot
\left(
\begin{array}{ccccccc}
2&2&1&1&2&2&1\\
3&3&2&2&4&3&3\\
5&4&3&3&6&5&4
\end{array}
\right)
\\
T_3^5 (p_{3s})&=&1_3^{\ot 13}\ot
\left(
\begin{array}{cccccccc}
2&2&1&1&1&2&2&1\\
3&3&2&2&2&3&3&3\\
5&4&3&4&3&6&5&4
\end{array}
\right)
\\
T_3^6 (p_{3s})&=&1_3^{\ot 15}\ot
\left(
\begin{array}{ccccccccc}
2&1&2&1&1&1&2&2&1\\
3&2&3&2&2&2&3&3&3\\
5&3&4&4&3&3&6&5&4
\end{array}
\right)
\end{eqnarray*}
With our notation this 3 body scattering is expressed as
\begin{eqnarray*}
&&
\zeta^0
\left[
\begin{array}{ccc}
2&2&2\\
3&3&3\\
4&4&5
\end{array}
\right]
\ot\zeta^6
\left[
\begin{array}{cc}
1&2\\
2&3\\
4&6
\end{array}
\right]
\ot\zeta^{11}
\left[
\begin{array}{c}
1\\
3\\
5
\end{array}
\right]
\\
\longrightarrow
&&
\zeta^9
\left[
\begin{array}{c}
2\\
3\\
5
\end{array}
\right]
\ot\zeta^5
\left[
\begin{array}{cc}
1&2\\
2&3\\
4&4
\end{array}
\right]
\ot\zeta^3
\left[
\begin{array}{ccc}
1&2&2\\
3&3&3\\
4&5&6
\end{array}
\right].
\end{eqnarray*}
Let us check that the RHS is obtained by applying either 
$(\tilde{R} \ot 1)(1\ot \tilde{R})(\tilde{R} \ot 1)$ or
$(1\ot \tilde{R})(\tilde{R} \ot 1)(1\ot \tilde{R})$ on the LHS.

\begin{eqnarray*}
LHS&\stackrel{\tilde{R}\ot1}{\longrightarrow}&
\zeta^{6-3}
\left[
\begin{array}{cc}
2&2\\
3&3\\
4&5
\end{array}
\right]
\ot\zeta^{0+3}
\left[
\begin{array}{ccc}
1&2&2\\
2&3&3\\
4&4&6
\end{array}
\right]
\ot\zeta^{11}
\left[
\begin{array}{c}
1\\
3\\
5
\end{array}
\right]\\ 
&\stackrel{1\ot\tilde{R}}{\longrightarrow}&
\zeta^3
\left[
\begin{array}{cc}
2&2\\
3&3\\
4&5
\end{array}
\right]
\ot\zeta^{11-0}
\left[
\begin{array}{c}
1\\
2\\
4
\end{array}
\right]
\ot\zeta^{3+0}
\left[
\begin{array}{ccc}
1&2&2\\
3&3&3\\
4&5&6
\end{array}
\right]\\ 
&\stackrel{\tilde{R}\ot1}{\longrightarrow}&
\zeta^{11-2}
\left[
\begin{array}{c}
2\\
3\\
5
\end{array}
\right]
\ot\zeta^{3+2}
\left[
\begin{array}{cc}
1&2\\
2&3\\
4&4
\end{array}
\right]
\ot\zeta^3
\left[
\begin{array}{ccc}
1&2&2\\
3&3&3\\
4&5&6
\end{array}
\right]
\end{eqnarray*}

\begin{eqnarray*}
LHS&\stackrel{1\ot\tilde{R}}{\longrightarrow}&
\zeta^0
\left[
\begin{array}{ccc}
2&2&2\\
3&3&3\\
4&4&5
\end{array}
\right]
\ot\zeta^{11-0}
\left[
\begin{array}{cc}
1\\
2\\
4
\end{array}
\right]
\ot\zeta^{6+0}
\left[
\begin{array}{cc}
1&2\\
3&3\\
5&6
\end{array}
\right]\\
&\stackrel{\tilde{R}\ot1}{\longrightarrow}&
\zeta^{11-2}
\left[
\begin{array}{c}
2\\
3\\
5
\end{array}
\right]
\ot\zeta^{0+2}
\left[
\begin{array}{ccc}
1&2&2\\
2&3&3\\
4&4&4
\end{array}
\right]
\ot\zeta^6
\left[
\begin{array}{cc}
1&2\\
3&3\\
5&6
\end{array}
\right]\\ 
&\stackrel{1\ot\tilde{R}}{\longrightarrow}&
\zeta^9
\left[
\begin{array}{c}
2\\
3\\
5
\end{array}
\right]
\ot\zeta^{6-1}
\left[
\begin{array}{cc}
1&2\\
2&3\\
4&4
\end{array}
\right]
\ot\zeta^{2+1}
\left[
\begin{array}{ccc}
1&2&2\\
3&3&3\\
4&5&6
\end{array}
\right]
\end{eqnarray*}
\end{example}

\subsection{Proof of the main theorem}

By Lemma \ref{lem:T_l commutes} it is sufficient to show the theorem for the highest weight 
element $p$ as a $U_q(sl_k)\times U_q(sl_{n-k})$-crystal, i.e., $e_ip=0$ for all $i\in I\setminus\{0,k\}$.
The following lemma is clear from the signature rule.

\begin{lemma}
The highest weight elements among two soliton states are given by 
\[
\tilde{p}_{2s}=[{\bf 0}]^{\ot c_1}\ot[{\bf 1}]^{\ot d_1}\ot[{\bf 0}]^{\ot(c_2-c_1-d_1)}
\ot[{\bf 3}]^{\ot \alpha}\ot[{\bf 2}_\pm] ^{\ot(d_2-\alpha-\beta)}\ot[{\bf 1}]^{\ot\beta}
\ot[{\bf 0}]^{\ot h_2},\nonumber
\]
where $h_2=L-c_2-d_2$ and $0\leq\alpha+\beta\leq d_2$.
\end{lemma}

By the commutation relation of time evolutions we get 
\[
T_r^t =T_{d_2 +1}^{-s} T_r^t T_{d_2 +1}^s \quad (t\gg s\gg 1,r>d_2).
\]
Thus it is sufficient to show the scattering rule by the time evolution $T_{d_2+1}$. 
We give the time evolution process for the case
$\alpha >\beta$ and $\alpha \leq \beta$ separately. Set $\delta =d_2 -\alpha +\beta$.

\begin{eqnarray*}
\lefteqn{{T_{d_2 +1}^t (\tilde{p}_{2s})} \quad (\alpha >\beta)} \\
&=&[{\bf 0}]^{\ot(c_1 +(d_2 +1)t)}\ot[{\bf 1}]^{\ot d_1}\ot[{\bf 0}]^{\ot(c_2 -c_1 -d_1-t)}
\ot[{\bf 3}]^{\ot\alpha}\ot[{\bf 2}_\pm]^{\ot(d_2 -\alpha -\beta)}\\
&&\ot[{\bf 1}]^{\ot\beta}\ot[{\bf 0}]^{\ot(h_2 -d_2 t)}\quad(0\leq t\leq c_2-c_1-d_1)\\
&=&[{\bf 0}]^{\ot(c_1 +(d_2 +1)t)}\ot[{\bf 1}]^{\ot(d_1 -s)}\ot[{\bf 4}]^{\ot s}\ot[{\bf 3}]^{\ot(\alpha -s)}
\ot[{\bf 2}_\pm]^{\ot(d_2 -\alpha -\beta)}\\
&&\ot[{\bf 1}]^{\ot\beta}\ot[{\bf 0}]^{\ot(h_2 -d_2 t)}\quad(0\leq s\leq\alpha -\beta,s=t-(c_2-c_1-d_1))\\
&=&[{\bf 0}]^{\ot(c_1 +(d_2 +1)t)}\ot[{\bf 1}]^{\ot(\delta -s)}\ot[{\bf 4}]^{\ot(\alpha -\beta)}
\ot[{\bf 3}]^{\ot\beta}\ot[{\bf 2}_\pm]^{\ot(d_2 -\alpha -\beta)}[{\bf 1}]^{\ot(\beta +s)}\\
&&\ot[{\bf 0}]^{\ot(h_1 -\delta -(d_2 +1)t)}\quad(0\leq s\leq d_1 -d_2,s=t-(c_2 -c_1 +d_2 -d_1 -\delta))\\
&=&[{\bf 0}]^{\ot(c_2 -\delta +d_2 t)}\ot[{\bf 1}]^{\ot(\delta +s)}\ot[{\bf 4}]^{\ot(\alpha -\beta -s)}
\ot[{\bf 3}]^{\ot(\beta +s)}\ot[{\bf 2}_\pm]^{\ot(d_2 -\alpha -\beta)}\\
&&\ot[{\bf 1}]^{\ot(\beta +d_1 -d_2)}\ot[{\bf 0}]^{\ot(h_1 -\delta -(d_2 +1)t)}
\quad(0\leq s\leq\alpha -\beta,s=t-(c_2 -c_1 -\delta))\\
&=&[{\bf 0}]^{\ot(c_2 -\delta +d_2 t)}\ot[{\bf 1}]^{\ot d_2}\ot[{\bf 0}]^{\ot s}\ot[{\bf 3}]^{\ot\alpha}
\ot[{\bf 2}_\pm]^{\ot(d_2 -\alpha -\beta)}\ot[{\bf 1}]^{\ot(d_1 -d_2 +\beta)}\\
&&\ot[{\bf 0}]^{\ot(h_1 -\delta -(d_2 +1)t)}\quad(0\leq s,s=t-(c_2 -c_1 +d_2 -2\delta))\\
\lefteqn{T_{d_2 +1}^t (\tilde{p}_{2s})\quad(\alpha\leq\beta)} \\
&=&[{\bf 0}]^{\ot(c_1 +(d_1 +1)t)}\ot[{\bf 1}]^{\ot d_1}\ot[{\bf 0}]^{\ot(c_2 -c_1 -d_1 -t)}
\ot[{\bf 3}]^{\ot\alpha}\ot[{\bf 2}_\pm]^{\ot(d_2 -\alpha -\beta)}\\
&&\ot[{\bf 1}]^{\ot\beta}\ot[{\bf 0}]^{\ot(h_2 -d_2 t)}\quad(0\leq t\leq c_2 -c_1 +d_2 -d_1 -\delta)\\
&=&[{\bf 0}]^{\ot(c_1 +(d_2 +1)t)}\ot[{\bf 1}]^{\ot(d_1 -s)}\ot[{\bf 0}]^{\ot(\beta -\alpha)}
\ot[{\bf 3}]^{\ot\alpha}\ot[{\bf 2}_{\pm}]^{\ot(d_2 -\alpha -\beta)}\ot[{\bf 1}]^{\ot(\beta +s)}\\
&&\ot[{\bf 0}]^{\ot(h_1 -\delta -(d_2 +1)t)}\quad(0\leq s\leq d_1 -d_2,s=t-(c_2 -c_1 -\delta))\\
&=&[{\bf 0}]^{\ot(c_2 -\delta +d_2 t)}\ot[{\bf 1}]^{\ot d_2}
\ot[{\bf 0}]^{\ot(\beta -\alpha +s)}\ot[{\bf 3}]^{\ot\alpha}\ot[{\bf 2}_\pm]^{\ot(d_2 -\alpha -\beta)}\\
&&\ot[{\bf 1}]^{\ot(\beta +d_1 -d_2)}\ot[{\bf 0}]^{\ot(h_1 -\delta -(d_2 +1)t)}
\quad(0\leq s,s=t-(c_2 -c_1 -\delta))
\end{eqnarray*}

{}From these we have 
\[
\zeta^{c_1} [{\bf 1}^{d_1}]\ot\zeta^{c_2} [{\bf 1}^{\beta} {\bf 2}_{\pm}^{d_2 -\alpha -\beta} {\bf 3}^{\alpha}]
\rightarrow \zeta^{c_2 -\delta} [{\bf 1}^{d_2}] \ot \zeta^{c_1 +\delta} [{\bf 1}^{d_1 -d_2 +\beta}
{\bf 2}_{\pm}^{d_2 -\alpha -\beta} {\bf 3}^{\alpha}].
\]
Recall the notation introduced in (\ref{eq:u split}).
{}From Theorem \ref{th:comb R} and Definition \ref{def:H}, we get the following.

\begin{lemma}
\begin{eqnarray*}
&&R([{\bf 1}^{d_1}]_k\ot [{\bf 1}^{\beta} {\bf 2}_{\pm}^{d_2 -\alpha -\beta} {\bf 3}^{\alpha}]_k)
=[{\bf 1}^{d_2}]_k \ot [{\bf 1}^{d_1 -d_2 +\beta}{\bf 2}_{\pm}^{d_2 -\alpha -\beta} {\bf 3}^{\alpha}]_k, \\
&&H([{\bf 1}^{d_1}]_k \ot [{\bf 1}^{\beta} {\bf 2}_{\pm}^{d_2 -\alpha -\beta} {\bf 3}^{\alpha}]_k)
=(-\frac{d_2 -\beta}{2} -\frac{\alpha}{2}) \mp (-\frac{d_2 -\beta}{2} +\frac{\alpha}{2}), \\
&&R([{\bf 1}^{d_1}]_{< k}\ot [{\bf 1}^{\beta} {\bf 2}_{\pm}^{d_2 -\alpha -\beta} {\bf 3}^{\alpha}]_{< k})
=[{\bf 1}^{d_2}]_{< k}\ot [{\bf 1}^{d_1 -d_2 +\beta}{\bf 2}_{\pm}^{d_2 -\alpha -\beta} {\bf 3}^{\alpha}]_{< k}, \\
&&H([{\bf 1}^{d_1}]_{< k}\ot [{\bf 1}^{\beta} {\bf 2}_{\pm}^{d_2 -\alpha -\beta} {\bf 3}^{\alpha}]_{< k})
=(-\frac{d_2 -\beta}{2} -\frac{\alpha}{2}) \pm (-\frac{d_2 -\beta}{2} +\frac{\alpha}{2}).
\end{eqnarray*}
\end{lemma}
{}From this lemma we get the phase shift as follows:
\[
2d_2 +H([{\bf 1}^{d_1}]_k \ot [{\bf 1}^{\beta} {\bf 2}_{\pm}^{d_2 -\alpha -\beta} {\bf 3}^{\alpha}]_k)
+H([{\bf 1}^{d_1}]_{< k} \ot [{\bf 1}^{\beta} {\bf 2}_{\pm}^{d_2 -\alpha -\beta} {\bf 3}^{\alpha}]_{< k})
=d_2 -\alpha +\beta =\delta,
\]
and see that have the change of the internal degree of freedom is given by the combinatorial $R$'s.

\section{Discussion}
In this paper we defined one soliton state $p\in B_1^{\ot L}$ so that it satisfies $E_1 (p)=1$. 
See Definition \ref{def:one soliton} and Remark \ref{rem:E_1=1}. But in the $E_1 (p)>1$ case also, 
there exist states in which a soliton-like particle moves without changing the internal degree of
freedom. As such examples we have 
\[
(a)\;
\left[
\begin{array}{cc}
2&1\\
4&6
\end{array}
\right]\quad
(b)\;
\left[
\begin{array}{ccc}
2&3&2\\
4&5&5
\end{array}
\right]\quad
(c)\;
\left[
\begin{array}{ccccc}
2&1&2&3&2\\
4&3&4&5&5
\end{array}
\right]
\]
How do we regard these states? Isn't it a one soliton state? 
Our answer is no. 
For a state $p$ define the numbers $N_d (p)$ ($d=1,2,\ldots$) by solving 
\[
E_l(p)=\sum_{d\geq 1} \min(d,l) N_d(p)\quad(l=1,2,\ldots),
\]
where $E_l(p)$ was defined in (\ref{eq:def_E}). 
Then $N_d$ gives the number of solitons of length $d$ for an $m$ soliton state
\[
\ldots[d_1]\ldots[d_2]\ldots \cdots \ldots[d_m]\ldots .
\]
if any adjacent solitons are separated enough. For the above examples the numbers $N_d$ reads as
\[
(a)\;N_1=2\quad(b)\;N_2=2\quad(c)\;N_2=3
\]
and the others are 0. Thus one may consider them as composite states of several solitons. 
However, it is in general difficult to identify the internal degree of freedom of each soliton,
and leave it as a future problem.

\vskip0.5cm\noindent
{\bf Acknowledgments}
\smallskip\par\noindent
The author would like to thank E. Date, K. Fukuda, M. Noumi, and H. Yamane for helpful discussions. 
Special thanks are due to M. Okado for his kind and warm guidance.

\end{document}